\documentclass[a4paper, 12pt, times]{article}
\usepackage{moreverb}
\usepackage{xcolor}
\usepackage{soul, authblk}
\usepackage{geometry} 
\usepackage[colorlinks,bookmarksopen,citecolor=red,urlcolor=red]{hyperref}
\usepackage{setspace}
\usepackage{amsmath}
\usepackage{nomencl}

\usepackage{natbib}
\usepackage[inkscapeformat=png]{svg}
\date{May 14, 2026}
\begin{document}


\title{Automatic differentiation in finite element stress updating}

\author{Mao Ouyang, Alexandros Petalas,
	William M. Coombs, Charles E. Augarde}

\affil{Department of Engineering, Durham University, South Road, Durham, DH1 3LE, County Durham, UK}

\maketitle
\begin{abstract}
	Computational solid mechanics relies on the discretisation of differential equations, requiring numerical differentiation and integration, often via quadrature. Recently, automatic differentiation (AD) has attracted interest as a tool to support these operations. This paper introduces the use of AD in nonlinear finite element analysis, focusing on stress updating, a key and repeatedly executed step in which stresses are computed from strain increments. Both implicit and explicit schemes require derivatives of constitutive models, which can be complex, error‑prone, and time‑consuming to derive analytically. We compare traditional analytical differentiation with automatic differentiation for implementing a hyperplastic Critical State model widely used in geotechnical engineering. Drained triaxial compression tests with varying overconsolidation ratios are simulated using both approaches. The results show that AD significantly simplifies the implementation of backward Euler stress integration by removing the need for manual derivation of derivatives. Despite this simplification, AD maintains high accuracy and robustness comparable to analytical approaches. The findings demonstrate that automatic differentiation can streamline the development of nonlinear material models, enabling efficient and reliable implementation of constitutive laws of arbitrary complexity. This opens the way for more flexible and maintainable computational mechanics codes.
\end{abstract}


\section{Introduction}
Finite element (FE) methods are widely used and developed by researchers in computational solid mechanics, and if modelling material non-linearity then a key aspect of any FE code is the determination of the updated stress at a Gauss point following a strain increment and its division into elastic and plastic parts. This stress updating requires differential terms linked to the constitutive model used which, if complex, can be challenging to obtain.

There are three ways to obtain derivatives: numeric, using finite differences which suffers from round-off errors; symbolic, where the derivative is formed as if we were doing it by hand and which has significant memory needs; and automatic which can provide the exact derivatives for complex functions by a simple, succinct, and clear approach \citep{vigliotti2021automatic}, which is fundamentally based on use of the chain rule. Automatic differentiation (AD)  was originally developed in the 1970s but has become of intense interest in machine learning in recent years \citep{margossian2019review}.
\nomenclature{AD}{Automatic Differentiation}
In solid mechanics, AD techniques have been implemented for hyperelastic constitutive models, e.g. to obtain the stress and tangent stiffness matrix from a strain energy function \citep{rothe2015automatic, vigliotti2021automatic}, recently for elasto-plasticity \citep{dummer2024robust, zhang2025si} and for finite strain generalised continuum models \citep{tanaka2016implementation}. 

This paper presents an introduction to AD for engineers developing FE codes, showing the ease with which it can be employed to update stress in a general (but fairly complex) geotechnical model as an example. The model used here is a hyperplastic Critical State model \citep{houlsby2006principles} with Willam-Warnke \citep{willam1974consitutive} Lode Angle Dependency (LAD), a combination proposed by \citet{coombs2011algorithmic}.
We demonstrate the accuracy that AD can achieve compared with using analytical derivatives in modelling Consolidated Drained (CD) triaxial compression tests with different overconsolidation ratios (OCRs). The stability and performance are also  investigated.
\nomenclature{CD}{Consolidated Drained}
\nomenclature{OCR}{Over-Consolidation Ratio}
\nomenclature{LAD}{Lode Angle Dependency}
A freely available package \textit{Enzyme} is used in this work for the AD calculations due to its high efficiency \citep{moses2020instead, moses2021reverse}.

\subsection{Constitutive model}\label{sec:model}
Here we use the hyperplastic Critical State model of  \cite{coombs2011algorithmic} which is based on a free energy function \citep{Houlsby1985} and a dissipation function \citep{Collins2002}. The elastic component of the free energy function defines a stress-elastic strain relationship with a variable bulk modulus, controlled by a reference elastic volumetric strain, $\varepsilon_{v0}^e$, at  reference pressure, $p_r$, and an elastic compressibility index, $\kappa$, that defines the slope of the elastic loading/unloading line in void ratio-hydrostatic pressure space, and a constant shear modulus, $G$. The dissipation function defines both the yield surface and the plastic flow rule, requiring five parameters: the slope of the Critical State Line in hydrostatic pressure-deviatoric stress ($p$-$q$) space, $M$, the ratio of the deviatoric yield radius under triaxial extension to that under triaxial compression, $\overline{\rho}_e$; two shape parameters that control the form of the yield surface in $p$-$q$ stress space, $\alpha$ and $\gamma$; and a history-dependent preconsolidation pressure that defines the hydrostatic extent of the the yield surface, $p_c$. The free energy function also includes a plastic component that controls the hardening/softening of the yield surface, which introduces a final parameter - the plastic compressibility index, $\lambda$, which defines the gradient of the virgin consolidation line in void ratio-hydrostatic pressure space. The yield surface shape parameters provide significant flexibility in the $p$-$q$ form of the yield function, whilst allowing recovery of the Modified Cam Clay (MCC) shape when $\alpha=\gamma=1$.  In total $9$ model constants are required\footnote{Note that the  reference elastic volumetric strain is often set to zero, $\varepsilon_{v0}^e=0$, and the reference pressure is usually taken as the sea-level atmospheric pressure, that is $p_r\approx 100$kPa. }, and one state variable, i.e., the preconsolidation pressure $p_c$, which evolves with plastic strain.

\nomenclature{BE}{Backward Euler}
\nomenclature{$\varepsilon_{v0}^e$}{Elastic volumetric strain at reference pressure}
\nomenclature{$p_r$}{Reference pressure}
\nomenclature{$\kappa$}{Elastic compressibility index}
\nomenclature{$G$}{Shear modulus}
\nomenclature{$p$}{Hydrostatic pressure}
\nomenclature{$q$}{Deviatoric stress}
\nomenclature{$M$}{Slope of Critical State line in the $p$-$q$ space}
\nomenclature{$\overline{\rho}_e$}{Ratio of the deviatoric radius under extension/compression}
\nomenclature{$\alpha,\ \gamma$}{Yield surface $p$-$q$ shape parameters}
\nomenclature{$p_c$}{Preconsolidation pressure; hydrostatic extent of the yield surface}
\nomenclature{$\lambda$}{Plastic compressibility index; slope of the virgin consolidation line}

\cite{coombs2011algorithmic} implemented this Critical State model using an implicit elastic predictor, plastic corrector Backward Euler (BE) algorithm and demonstrated the accuracy and robustness of this approach even under large strain increments\footnote{For brevity, details of the stress integration algorithm are not provided here but they are identical to the procedure outlined in Fig. 11 of} \cite{coombs2011algorithmic}. Although efficient and robust, there is considerable overhead in developing such an implementation due to the requirement to obtain first, and sometimes second, linearisations of the  residual equations with respect to the primary unknowns: the principal elastic strains, preconsolidation pressure and plastic multiplier in the case of \cite{coombs2011algorithmic}. 
This tedious process can be avoided entirely by using AD in the BE algorithm. This still requires the definition of residual equations (which specify when a valid state has been obtained) and associated primary unknowns, but the linearisation of these residuals with respect to the unknowns is now automatic. A key potential benefit of AD-enabled stress updating is it allows rapid exploration of variations in constitutive \emph{ingredients}, such as hardening laws, elasticity relationships, etc., without redetermining and reimplementing the linearised residual equations.  Below we compare the accuracy, stability and performance of this approach versus the original BE algorithm of \cite{coombs2011algorithmic}.      
\newpage
\section{Results}
In this section, our focus is on three aspects: (a) accuracy, (b) stability, and (c) performance. 
\begin{figure*}[htpb!]
	\centering
	\includegraphics[width=1\textwidth]{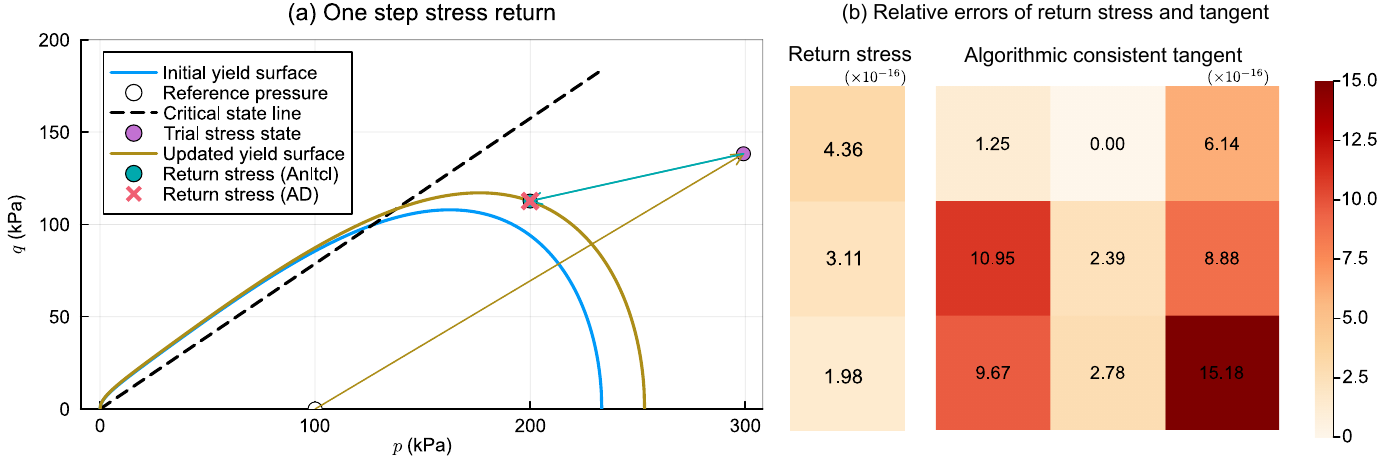}
	\caption{(a) One step stress return with analytical derivatives (\textit{Anltcl}) and automatic differentiation (\textit{AD});
		(b) Relative errors of return stress and algorithmic consistent tangent;
		The relative errors are calculated by Equation \eqref{eq:e}.}
	\label{fig:accuracy}
\end{figure*}
\subsection{Accuracy}
The accuracy of the two algorithms is assessed by comparing the final stress solution after applying a strain increment in a single step, which activates the yield surface (YS) and is subsequently corrected back onto it. The principal strain increment components are given by \{$d\boldsymbol{\varepsilon}$\}= \{16, 0, -8\} $\times$ $10^{-3}$. The material constants for the hyperplastic model of \citet{coombs2011algorithmic} are as follows: $G$=4000 kPa, $\kappa$=0.0073, $\varepsilon_{v0}^e$=0, $p_r$=100 kPa, $\lambda$=0.0447, $M$=0.964, $\overline{\rho}_e$=0.729, $\alpha$=0.3, $\gamma$=0.9, $p_c$=233.3 kPa. The initial stress state, prior to the application of the strain increment, is given by $p$=100 kPa, $q$=0 kPa.

\nomenclature{YS}{Yield Surface}

Figure \ref{fig:accuracy}(a) shows the calculated stress for the applied strain increment for both integration algorithms. Figure \ref{fig:accuracy}(b) quantifies the error of the algorithm using AD, assuming analytical derivatives to give the ``correct" solution. Figure \ref{fig:accuracy}(a) shows that both algorithms, after the ``Trial" elastic prediction (purple dot), predict an almost identical ``Return" stress (green dot and red cross). This indicates that the algorithm does not lose accuracy when AD is used compared to the traditional process with analytical derivatives.  

More specifically, the relative \textit{error} is quantified for each principal component of the predicted stress tensor and for each component of the estimated algorithmic consistent tangent matrix \footnote{The  algorithmic consistent tangent is the consistent linearisation of the constitutive algorithm, which is used to construct the element stiffness in implicit finite element analysis, full details are given in }\cite{coombs2011algorithmic}. Specifically  the algorithmic tangent stiffness is calculated using Eqns (14) and (68) of that reference. The entire procedure is identical apart from the fact that derivatives are replaced by AD here, using,  
\begin{equation}
	\label{eq:e}
	\textit{error} = \frac{|\textit{AD} - \textit{Anltcl}|}{\textit{Anltcl}},
\end{equation}
where \textit{AD} and \textit{Anltcl} represent the results obtained using AD and analytical derivatives, respectively. The relative \textit{error} of the algorithm using AD for each stress component and for each component of the algorithmic tangent matrix is on the order of \(10^{-16}\), again confirming that using AD does not sacrifice accuracy.

To further examine the accuracy of the algorithm using AD, the two integration algorithms are used in a 3D Finite Element code \citep{coombs201070} for the simulation of isotropically consolidated, drained triaxial tests for different OCRs. An 8-noded hexahedral element with 8 Gauss points is used. The material parameters follow \citet{coombs2011algorithmic} and are the same as the previous section, except for $p_c$ which evolves with the plastic straining.

All four nodes of the lower horizontal face are constrained for vertical displacement, while the four nodes on two neighbouring vertical faces are constrained for horizontal displacements in their normal directions. On the remaining two vertical faces, a uniform hydrostatic pressure is applied during the consolidation stage, where an initial reference pressure of \(p_r = 100\) kPa is applied to represent the sea-level atmospheric pressure. The cube is then consolidated to \(p = 233.3\) kPa in 40 load steps. 

After reaching the consolidation stress \( p = 233.3 \) kPa, the cube is unloaded by isotropic expansion to \( p = 233.3, 186.7, 155.6, 166.7, 58.3, \) and \( 23.3 \) kPa, corresponding to OCRs of 1, 1.25, 1.5, 2, 4, and 10, respectively, in 40 load steps\footnote{Note that for OCR=1, the boundary conditions are kept constant. We use 40 load steps here for the general purpose and comprehensive comparisons.}. Then, an axial vertical displacement is applied to the  top horizontal face, with a final applied axial strain of \( \varepsilon_a = 0.14 \).  

\begin{figure*}[htpb!]
	\centering
	\includegraphics[width=1\textwidth]{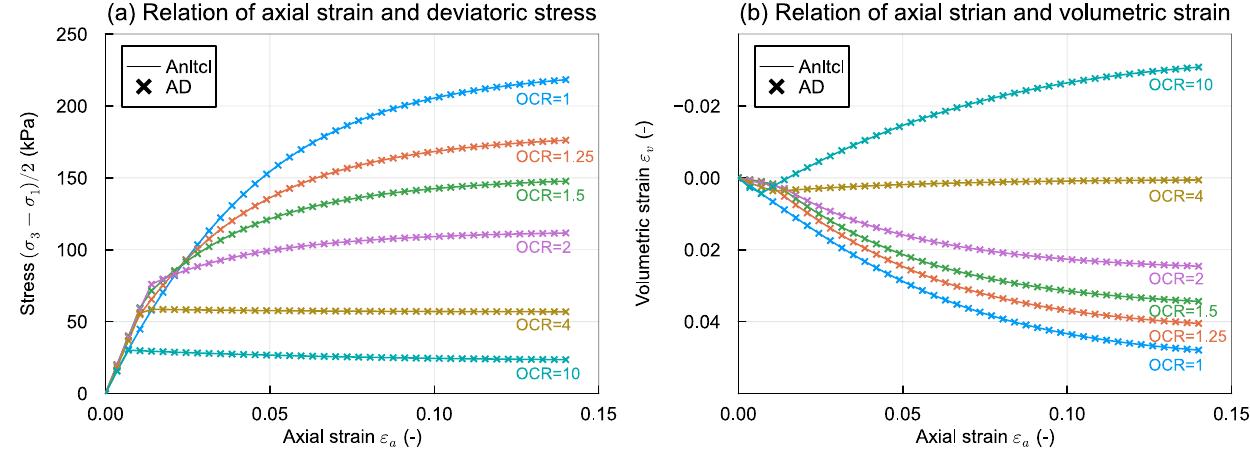}
	\caption{Simulation results, using the implicit Backward Euler stress integration algorithm with analytical derivatives (solid lines) and automatic differentiation (cross symbols). (a) Deviatoric stress versus axial strain;
		(b) Volumetric strain versus axial strain.}
	\label{fig:applicability}
\end{figure*}

Reference results are obtained by applying the traditional BE algorithm with analytical derivatives \citep{coombs2011algorithmic}, shown as solid lines in Fig. \ref{fig:applicability}. The symbols (crosses) represent the  solutions with AD. The results demonstrate that  AD provides an identical (i.e. within machine precision) solution to the reference algorithm, and thus we conclude that the accuracy when using AD, within a finite element, is also not compromised. 
\nomenclature{$\varepsilon_a$}{Axial strain}

\subsection{Stability}
\begin{figure*}[htpb!]
	\centering
	\includegraphics[width=1\textwidth]{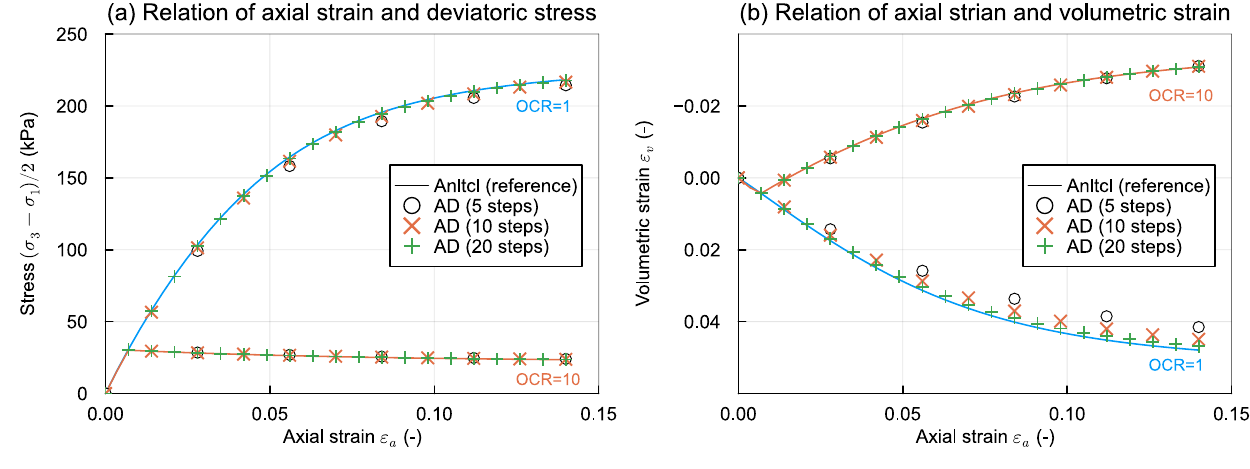}
	\caption{Stability of the code with AD by different load steps. The solid lines are the reference solutions with the analytical derivatives. The scatters denote the results with AD by 5, 10 and 20 load steps.}
	\label{fig:stability}
\end{figure*}

To investigate the stability and robustness of the BE algorithm with AD, we perform two  simulations of the previous section, choosing two OCRs (1 and 10) with different numbers of steps, i.e., 5, 10, and 20. The boundary conditions and the applied vertical displacement remain the same as the simulations above.  The results in Fig. \ref{fig:stability} suggest that the algorithm with AD remains stable even with a very small number of steps (5 steps), with an acceptably small loss of accuracy.
\subsection{Performance}
\begin{figure*}[htpb!]
	\centering
	\includegraphics[width=1.00\textwidth]{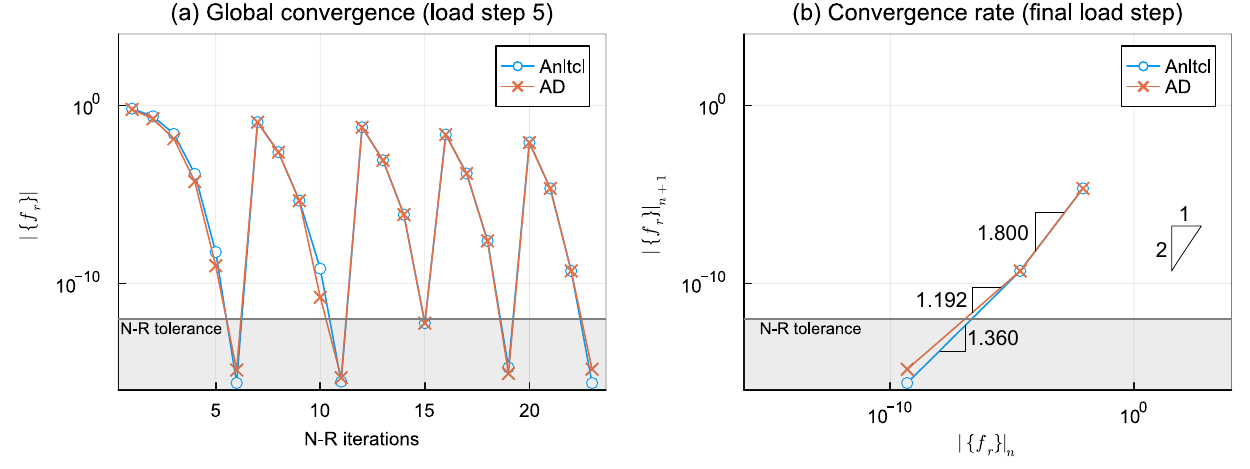}
	\caption{Convergence results of the finite element code with AD for the case with 5 load steps.
		(a) Norm of the residual out-of-balance force for the drained compression tests;
		(b) Norm of the residual out-of-balance force against the previous out-of-balance force for the final load step.}
	\label{fig:convergence}
\end{figure*}

The performance of the BE algorithm with AD is evaluated by assessing the local convergence of the Newton-Raphson (N-R) iterative procedure during stress correction at a Gauss point, as well as the global convergence when using a consistent tangent matrix computed with AD. The performance is again compared to the BE algorithm using  analytical derivatives.
\nomenclature{N-R}{Newton-Raphson}
The convergence of the global BE iterations is shown in Fig.~\ref{fig:convergence}(a) for the simulation with $\text{OCR}=1$ and 5 load steps. The normalised out-of-balance force, $|\{f_r\}|$, is defined as
\begin{equation} 
|\{f_r\}| = \frac{\sqrt{(\{f_{ext}\} - \{f_{int}\})^T(\{f_{ext}\} - \{f_{int}\}})}{\sqrt{\{f_{ext}\}^T\{f_{ext}\}}} 
\end{equation}
where ${f_{ext}}$ and ${f_{int}}$ are the external and internal force vectors, respectively. The convergence tolerance is set to $10^{-12}$. The convergence of the global iterative N-R scheme, when using a consistent tangent matrix calculated by AD, is not distinguishably affected. As shown in Fig.\ref{fig:convergence}(a), the out-of-balance force decreases at a similar rate. More specifically, in Fig.\ref{fig:convergence}(b), the norm of the residual out-of-balance force at the current iteration is plotted against the corresponding value from the previous iteration. In both cases, the initial convergence rate is 1.8 (superlinear convergence). In the final iteration, the case with the analytically consistent tangent matrix gives a convergence rate of 1.360, whereas the case with AD gives a convergence rate of 1.192 - in both cases the rate is limited by double precision accuracy.

\nomenclature{$f_r$}{Out-of-balance force}
\nomenclature{$f_{ext}$}{External force}
\nomenclature{$f_{int}$}{Internal force}

\begin{figure}[htpb!]
	\centering
	\includegraphics[width=0.48\textwidth]{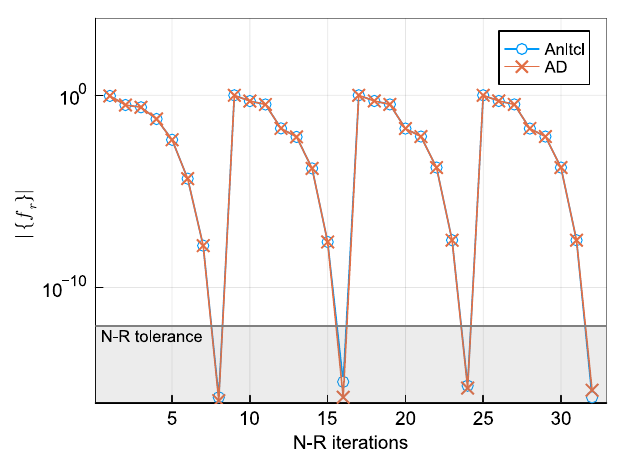}
	\caption{Local convergence at one Gauss point of the final load step.}
	\label{fig:localconvergence}
\end{figure}

The local convergence behaviour for the two cases is presented in Fig.~\ref{fig:localconvergence}, focusing on the last four global iterations of Step 5. Each global iteration triggers a series of local iterations, performed by applying a strain increment to the constitutive equations, which are integrated using the BE algorithm. As shown in Fig.~\ref{fig:localconvergence}, the convergence behaviour of the algorithm employing AD is nearly identical to that obtained using analytical derivatives. It can therefore be concluded that both local and global convergence are not notably affected by the use of AD.
A final aspect of performance is is computational efficiency, and in this study we have found that the analytical derivative approach is slightly quicker than the AD approach for the examples present. However, this should not demotivate its use. 
\citet{griewank2008evaluating} noted that in terms of computational complexity, AD guarantees that the amount of arithmetic increases by no more than a small constant factor and with the rapid development of AD techniques, e.g., the Low-Level Virtue Machine (LLVM), and hardware, e.g., Graphics Processing Units (GPU) and Tensor Processing Units (TPU), it is likely this will cease to be a problem in future.
\nomenclature{GPU}{Graphic Processing Units}
\nomenclature{TPU}{Tensor Processing Units}
\nomenclature{LLVM}{Low-Level Virtue Machine}
\section{Summary}
Non-linear finite element analysis in solid mechanics involves many operations to determine derivatives of functions,  for updating stresses and tangent stiffness matrices.
AD offers a simple way to explore and implement complex constitutive models 
without affecting accuracy or stability. In this paper we have focussed on those beloved of geotechnical engineers and while many constitutive models have smooth equations, multiple state variables with non-linear evolution laws make analytical consistent tangents for implicit integration practically infeasible. Our methodology applies directly to such cases, and the model used here is sufficiently complex to demonstrate this. Established yield-surface smoothing procedures \cite{ABBO1995, COOMBS2016, ALCAZAR2026} also work seamlessly with AD. Moreover, AD has proven effective in dealing even with functions with discontinuities, further supporting its suitability for elasto-plasticity.   

\setlength{\nomitemsep}{-\parsep}

\bibliography{automaticdifferentiation}

@article{COOMBS2016,
title = {{NURBS} plasticity: Yield surface representation and implicit stress integration for isotropic inelasticity},
journal = {Computer Methods in Applied Mechanics and Engineering},
volume = {304},
pages = {342-358},
year = {2016},
author = {William M. Coombs and Oscar A. Petit and Yousef {Ghaffari Motlagh}},
}

@article{ALCAZAR2026,
title = {Topology optimization for pressure-dependent elastoplastic structures considering a smooth hyperbolic approximation of the {D}rucker-{P}rager yield criterion},
journal = {Computer Methods in Applied Mechanics and Engineering},
volume = {449},
pages = {118446},
year = {2026},
author = {Emily Alcazar and Glaucio H. Paulino}}

@article{ABBO1995,
title = {A smooth hyperbolic approximation to the {M}ohr-{C}oulomb yield criterion},
journal = {Computers \& Structures},
volume = {54},
number = {3},
pages = {427-441},
year = {1995},
author = {A.J. Abbo and S.W. Sloan}
}

@Article{zhang2025si,
  author  = {Zhang, Pin and Yin, Zhen-Yu and Sheil, Brian},
  journal = {Géotechnique Letters},
  title   = {{si-PiNet: a novel stress integration method for elastoplastic models}},
  year    = {2025},
  number  = {1},
  pages   = {12--18},
  volume  = {15},
}

@article{Houlsby1985,
title = {The use of a variable shear modulus in elastic-plastic models for clays},
journal = {Computers \& Geotechnics},
volume = {1},
number = {1},
pages = {3-13},
year = {1985},
author = {G.T. Houlsby}
}

@article{Collins2002,
author = {Collins, Ian F. and Hilder, Tamsyn},
title = {A theoretical framework for constructing elastic/plastic constitutive models of triaxial tests},
journal = {International Journal for Numerical and Analytical Methods in Geomechanics},
volume = {26},
number = {13},
pages = {1313-1347},
year = {2002}
}

@Article{coombs2011algorithmic,
  author  = {William M. Coombs and Roger S. Crouch},
  journal = {Computer Methods in Applied Mechanics and Engineering},
  title   = {Algorithmic issues for three-invariant hyperplastic {Critical State} models},
  year    = {2011},
  number  = {25},
  pages   = {2297--2318},
  volume  = {200},
}

@InProceedings{willam1974consitutive,
  author    = {Willam, K.~J. and Warnke, E.~P.},
  booktitle = {International Association of Bridge and Structural Engineers Seminar on Concrete Structures Subjected to Triaxial Stresses},
  title     = {Consitutive model for the triaxial behvaviour of concrete},
  year      = {1974},
}

@Article{vigliotti2021automatic,
  author  = {Vigliotti, Andrea and Auricchio, Ferdinando},
  journal = {Archives of Computational Methods in Engineering},
  title   = {Automatic differentiation for solid mechanics},
  year    = {2021},
  issn    = {1886-1784},
  number  = {3},
  pages   = {875--895},
  volume  = {28}
}

@Article{dummer2024robust,
  author  = {Alexander Dummer and Matthias Neuner and Peter Gamnitzer and Günter Hofstetter},
  journal = {Computer Methods in Applied Mechanics and Engineering},
  title   = {Robust and efficient implementation of finite strain generalized continuum models for material failure: Analytical, numerical, and automatic differentiation with hyper-dual numbers},
  year    = {2024},
  issn    = {0045-7825},
  pages   = {116987},
  volume  = {426},
}

@Article{rothe2015automatic,
  author  = {Rothe, Steffen and Hartmann, Stefan},
  journal = {Archive of Applied Mechanics},
  title   = {Automatic differentiation for stress and consistent tangent computation},
  year    = {2015},
  issn    = {1432-0681},
  number  = {8},
  pages   = {1103--1125},
  volume  = {85},
}

@InProceedings{moses2020instead,
  author    = {Moses, William and Churavy, Valentin},
  booktitle = {Advances in Neural Information Processing Systems},
  title     = {Instead of rewriting foreign code for machine learning, automatically synthesize fast gradients},
  year      = {2020},
  editor    = {H. Larochelle and M. Ranzato and R. Hadsell and M. F. Balcan and H. Lin},
  pages     = {12472--12485},
  volume    = {33},
}

@InProceedings{moses2021reverse,
  author    = {Moses, William S. and Churavy, Valentin and Paehler, Ludger and H\"{u}ckelheim, Jan and Narayanan, Sri Hari Krishna and Schanen, Michel and Doerfert, Johannes},
  booktitle = {Proceedings of the International Conference for High Performance Computing, Networking, Storage and Analysis},
  title     = {Reverse-mode automatic differentiation and optimization of {GPU} kernels via {Enzyme}},
  year      = {2021},
  address   = {New York, NY, USA},
  series    = {SC '21},
  articleno = {61},
  numpages  = {16},
}

@Article{tanaka2016implementation,
  author  = {Masato Tanaka and Daniel Balzani and Jörg Schröder},
  journal = {Computer Methods in Applied Mechanics and Engineering},
  title   = {Implementation of incremental variational formulations based on the numerical calculation of derivatives using hyper dual numbers},
  year    = {2016},
  issn    = {-},
  pages   = {216--241},
  volume  = {301},
}

@Book{houlsby2006principles,
  author    = {G. T. Houlsby and A. M. Puzrin},
  publisher = {Springer London},
  title     = {{Principles of Hyperplasticity}},
  year      = {2006},
}

@Article{margossian2019review,
  author  = {Margossian, Charles C.},
  journal = {WIREs Data Mining and Knowledge Discovery},
  title   = {A review of automatic differentiation and its efficient implementation},
  year    = {2019},
  number  = {4},
  pages   = {e1305},
  volume  = {9},
}

@Misc{coombs201070,
  author            = {Coombs, W.M. and Crouch, R.S. and Augarde, C.E.},
  title             = {{70-line 3D} finite deformation elastoplastic finite-element code},
  year              = {2010},
  conference        = {7th European Conference on Numerical Methods in Geotechnical Engineering (NUMGE)},
  editor            = {T., Benz and S., Nordal},
  isbn              = {9780415592390},
  pages             = {151--156},
  publicationstatus = {Published},
  publisher         = {Taylor and Francis},
}

@Book{griewank2008evaluating,
  author    = {Griewank, Andreas and Walther, Andrea},
  publisher = {Society for Industrial and Applied Mathematics},
  title     = {{Evaluating Derivatives}},
  year      = {2008},
  edition   = {Second},
}

\end{document}